\documentclass[11pt]{article}
\usepackage{fullpage}
\usepackage{amssymb,amsfonts}
\usepackage{amsthm} 
\usepackage[cmex10]{amsmath}
\usepackage[usenames,dvipsnames,svgnames]{xcolor} 
\usepackage{multirow}
\usepackage{graphicx}
\usepackage{booktabs}
\usepackage{url}

\usepackage[font={small,it}]{caption}

\usepackage{array} 
\usepackage{setspace}
\ifdefined\isdraft
\usepackage[sort,compress]{natbib}
\else
\usepackage{cite} 
\fi

\usepackage[utf8]{inputenc}
\usepackage[T1]{fontenc}
\usepackage{lmodern} 

\usepackage{verbatim}

\usepackage{hyperref}
\hypersetup{pdfauthor={Volkan Cevher, Stephen Becker, Mark Schmidt},
pdftitle={Convex Optimization for Big Data},
    colorlinks=true,
    citecolor=Mulberry,
    urlcolor=Bittersweet,
}
\usepackage[normalem]{ulem}

\DeclareMathOperator{\prox}{prox}

\DeclareMathOperator*{\argmin}{argmin}

\newcommand{\eps}{\varepsilon}

\newcommand{\R}{\mathbb{R}}

\newcommand{\E}{\operatorname{\mathbb{E}}}

\newcommand{\order}{\mathcal{O}}

\DeclareMathAlphabet{\mathbbb}{U}{bbold}{m}{n}

\newcommand{\defeq}{\stackrel{\text{\tiny def}}{=}}

\newcommand{\norm}[1]{{\left\lVert{#1}\right\rVert}}

\usepackage{algorithm,algorithmic}

\usepackage{chngcntr}
\counterwithin*{paragraph}{subsection} 

\newcommand{\iterate}[2]{#1^{#2}}  

\newcommand{\xk}{\iterate{x}{k}}
\newcommand{\xkk}{\iterate{x}{k+1}}
\newcommand{\zk}{\iterate{z}{k}}
\newcommand{\zkk}{\iterate{z}{k+1}}
\newcommand{\vk}{\iterate{v}{k}}
\newcommand{\vkk}{\iterate{v}{k+1}}
\newcommand{\uk}{\iterate{u}{k}}
\newcommand{\ukk}{\iterate{u}{k+1}}
\newcommand{\xzero}{\iterate{x}{0}}
\newcommand{\xStar}{\iterate{x}{\star}}
\newcommand{\vzero}{\iterate{v}{0}}
\newcommand{\zzero}{\iterate{z}{0}}
\newcommand{\uzero}{\iterate{u}{0}}

\def\Real{\mathbb{R}}

\title{Convex Optimization for Big Data}
\author{Volkan Cevher\footnote{\'Ecole Polytechnique F\'ed\'erale de Lausanne (EPFL)---Laboratory for Information and Inference Systems (LIONS)}, Stephen Becker\footnote{University of Colorado---Applied Math}, and Mark Schmidt\footnote{University of British Columbia---Laboratory for Computational Intelligence}
}
\date{September 2014}
\begin{document}
\maketitle
 % spmag_abstract
%\vspace{-3mm} 
This article reviews recent advances in convex optimization algorithms for Big Data, which aim to reduce the computational, storage, and communications bottlenecks. We provide an overview of  this emerging field, describe contemporary approximation techniques like first-order methods and randomization for scalability, and survey the important role of parallel and distributed computation. The new Big Data algorithms are based on surprisingly simple principles and attain staggering accelerations even on classical problems.

%\vspace{-3mm}

\subsection*{Convex optimization in the wake of Big Data}
{\tiny \pdfbookmark[1]{Introduction}{refSection} }
Convexity in signal processing dates back to the dawn of the field, with problems like least-squares being ubiquitous across nearly all sub-areas. However, the importance of convex formulations and optimization has increased even more dramatically in the last decade due to the rise of new theory for structured sparsity and rank minimization, and successful statistical learning models like support vector machines.
These formulations are now employed in a wide variety of signal processing applications including compressive sensing, medical imaging, geophysics, and bioinformatics \cite{wainwright2014structured,chandrasekaran2012convex,nesterov2013first,CombettesPesquetChapter}. 
 
There are several important reasons for this explosion of interest, with two of the most obvious ones being the existence of efficient algorithms for computing globally optimal solutions and the ability to use convex geometry to prove useful properties about the solution \cite{wainwright2014structured,chandrasekaran2012convex}. A unified convex formulation also transfers useful knowledge across different disciplines,
such as sampling and computation,  % "such as on" change... IGNORE
that focus on different aspects of the same underlying mathematical problem~\cite{chandrasekaran2013computational}.

However, the renewed popularity of convex optimization places convex algorithms under tremendous pressure to accommodate increasingly large data sets and to solve problems in unprecedented dimensions. Internet, text, and imaging problems (among a myriad of other examples) no longer produce data sizes from megabytes to gigabytes, but rather from terabytes to exabytes. Despite the progress in parallel and distributed computing, the practical utility of classical algorithms like interior point methods may not go beyond discussing the theoretical tractability of the ensuing optimization problems \cite{nesterov2013first}.

In response, convex optimization is reinventing itself for Big Data where the data and parameter sizes of optimization problems are too large to process locally, and where even basic linear algebra routines like Cholesky decompositions and matrix-matrix or matrix-vector multiplications that algorithms take for granted are prohibitive. In stark contrast, convex  algorithms also no longer need to seek high-accuracy solutions since Big Data models are necessarily simple or inexact  \cite{bottou2007tradeoffs}.
%\vspace{-3mm}

\subsection*{The basics}
We describe the fundamentals of Big Data optimization via the following composite formulation 
\begin{equation}\label{eq: master}
F^* \defeq \min_{x} \left\{F(x)\defeq f(x) + g(x) : x \in \R^p \right\},
\end{equation}
where $f$ and $g$ are convex functions. We review efficient numerical methods to  obtain an optimal solution $\xStar$ of \eqref{eq: master} 
as well as required assumptions on $f$ and $g$.
 Such composite convex minimization problems naturally arise in signal processing when we estimate unknown parameters $x_0 \in \mathbb{R}^p$ from data $y \in \mathbb{R}^n$. In \emph{maximum a posteriori} estimation, for instance, we regularize a smooth data likelihood function as captured by $f$ typically with a non-smooth prior term $g$ that encodes parameter complexity \cite{wainwright2014structured}. 

A basic understanding of Big Data optimization algorithms for \eqref{eq: master} rests on three key pillars:
\begin{itemize}
\item {\bf First-order methods} (Section I): First order methods obtain low- or medium-accuracy numerical solutions by using only first-order oracle information from the objective, such as gradient estimates. They can also handle the important \emph{non-smooth} variants of \eqref{eq: master} by making use of the \emph{proximal mapping} principle. These methods feature nearly dimension-independent convergence rates, they are theoretically robust to the approximations of their oracles, and they typically rely on computational primitives that are ideal for distributed and parallel computation. 
\item {\bf Randomization} (Section II): Randomization techniques particularly stand out among many other approximation techniques to enhance the scalability of first-order methods since we can control 
their expected behavior. Key ideas include random partial updates of optimization variables, replacing the deterministic gradient and proximal calculations with cheap statistical estimators, and speeding up basic linear algebra routines via randomization. 
\item {\bf Parallel and distributed computation} (Section III): First-order methods naturally provide a flexible framework to distribute optimization tasks and perform computations in parallel. 
Surprisingly, we can further augment these methods with approximations for increasing levels of scalability, from idealized synchronous parallel algorithms with centralized communications to enormously-scalable \emph{asynchronous} algorithms with \emph{decentralized} communications. 
\end{itemize}

The three concepts above complement each other to offer surprising scalability benefits for Big Data optimization. For instance, \emph{randomized} first-order methods can exhibit significant acceleration over their deterministic counterparts since they can generate a good quality solution with high probability by inspecting only a negligibly small fraction of the data \cite{nesterov2013first}. Moreover, since the computational primitives of such methods are inherently approximate, we can often obtain near linear speed-ups with a large number of processors \cite{niu2011hogwild,richtarik2012parallel}, which is a difficult feat when exact computation is required.

 % spmag_motivationChallenges
%\vspace{-3mm}%
\subsection*{A motivation for first-order methods}
A main source of Big Data problems is the ubiquitous linear observation model in many disciplines:
\begin{equation}\label{eq: linear model}
y= \Phi x_0 + z,
\end{equation}
where $x_0$ is an unknown parameter, $\Phi\in \mathbb{R}^{n\times p}$ is a known matrix and $z \in \mathbb{R}^{n}$ encodes unknown perturbations or noise---modeled typically with zero-mean iid Gaussian entries with variance $\sigma^2$. Linear observations sometimes arise directly from the basic laws of physics as in magnetic resonance imaging and geophysics problems. Other times, \eqref{eq: linear model} is an approximate model for more complicated nonlinear phenomena as in recommender systems and phase retrieval applications. 

The linear model \eqref{eq: linear model} along with low-dimensional signal models on $x_0$, such as sparsity, low total-variation, and low-rankness,  has been an area of intense research activity in signal processing. Hence, it is instructive to first study the choice of convex formulations and their scalability implications  here. 
The classical convex formulation in this setting has always been the least squares (LS) estimator
\begin{equation}\label{eq: LS}
\widehat{x}_{\rm LS} = \argmin_{x\in\Real^p} \left\{ F(x) := \frac{1}{2}\norm{y - \Phi x}_2^2 \right\},
\end{equation}
which can be efficiently solved by Krylov subspace methods using only matrix-vector multiplications. 
An important variant to \eqref{eq: LS} is the $\ell_1$-regularized least absolute shrinkage and selection operator (LASSO), which features the composite form \eqref{eq: master} 
\begin{equation}
\label{eq: lasso}
\widehat{x}_{\rm LASSO} = \argmin_{x\in\Real^p} \left\{ F(x) :=  \frac{1}{2}\norm{y - \Phi x}_2^2 + \lambda\norm{x}_1\right\},
\end{equation}
where $\lambda$ controls the strength of the regularization.  Compared to the LS estimator, the LASSO estimator has the advantage of producing \emph{sparse} solutions (i.e.,  $\widehat{x}_{\rm LASSO}$ has mostly zero entries), but its numerical solution  is essentially harder since the regularizing term is non-smooth. 

It turns out that the $\ell_1$-regularization in sparse signal recovery with the linear model \eqref{eq: linear model} is indeed critical when we are data deficient (i.e., $n<p$). Otherwise, the LASSO formulation  imparts only a denoising effect to the solution when $n\ge p$.  Theoretical justifications of the LASSO \eqref{eq: lasso} over the LS estimation \eqref{eq: LS} come from statistical analysis and the convex geometry of \eqref{eq: lasso}, and readily apply to many other low-dimensional signal models and their associated composite formulations \cite{wainwright2014structured,chandrasekaran2012convex,nesterov2013first}. 

\begin{table}[!t] %[!t]\small
\centering
 \begin{tabular}{rrrrrrr}
\toprule
Dimension & \multicolumn{2}{c}{Time} & \multicolumn{2}{c}{Error $\|\widehat{x}-x_0\|^2/\sigma^2$} &
\multicolumn{2}{c}{Iterations} \\
\cmidrule(r){2-3} \cmidrule(r){4-5} \cmidrule(r){6-7}
 & SDPT3 & TFOCS & SDPT3 & TFOCS &SDPT3 &TFOCS\\
\midrule
%\midrule
    128 &      0.3 s &      0.3 s  &    1.2 &   1.2&      10 &      94 \\ 
    512 &      2.2 s &      0.3 s  &    2.3 &   2.3&      11 &     121 \\ 
   1024 &     16.0 s &      0.5 s  &    2.4 &   2.4&      12 &     157 \\ 
   2048 &    145.0 s &      0.7 s  &    2.8 &   2.8&      12 &     234 \\ 
   4096 &      N/A &      1.0 s  &    N/A &   3.3&     N/A &     281 \\ 
  16384 &      N/A &      2.9 s  &    N/A &   3.7&     N/A &     527 \\ 
 131072 &      N/A &     40.2 s  &    N/A &   4.4&     N/A &    1265 \\ 
1048576 &      N/A &    838.5 s  &    N/A &   5.1&     N/A &    3440 \\ 
\bottomrule
\end{tabular}
    \caption{A numerical comparison of the default first-order method implemented in TFOCS~\cite{TFOCS}
    %(\url{http://cvxr.com/tfocs})
        versus the interior point method SDPT3 implemented in CVX~\cite{CVX} for the LASSO problem \eqref{eq: lasso} with $\lambda= 2 \sigma\sqrt{2\log p}$.  In the linear observation model \eqref{eq: linear model}, the matrix $\Phi$ is a randomly sub-sampled DCT matrix with $n=p/2$, the signal $x_0$ has $s=p/25$ non-zero coefficients with norm $\|x_0\|_2^2\approx s$, and the noise $z$ has variance $\sigma^2=10^{-4}$. 
}%\vspace{-2mm}
\label{table:cvx_vs_tfocs}
\end{table}

Table~\ref{table:cvx_vs_tfocs} illustrates key hallmarks of the first-order methods with the LASSO problem against the classical interior point method: nearly dimension-independent convergence and the ability to exploit implicit linear operators (e.g.,\ the DCT transform). In contrast, interior point methods require much larger space and have near cubic dimension dependence due to the application of dense matrix-matrix multiplications or Cholesky decompositions in finding the Newton-like search directions. Surprisingly, the LASSO formulation possesses additional structures that provably enhance the convergence of the first-order methods \cite{wainwright2014structured}, making them competitive in accuracy even to the interior point method.

\begin{figure}\centering
    \includegraphics[width=0.35 \textwidth]{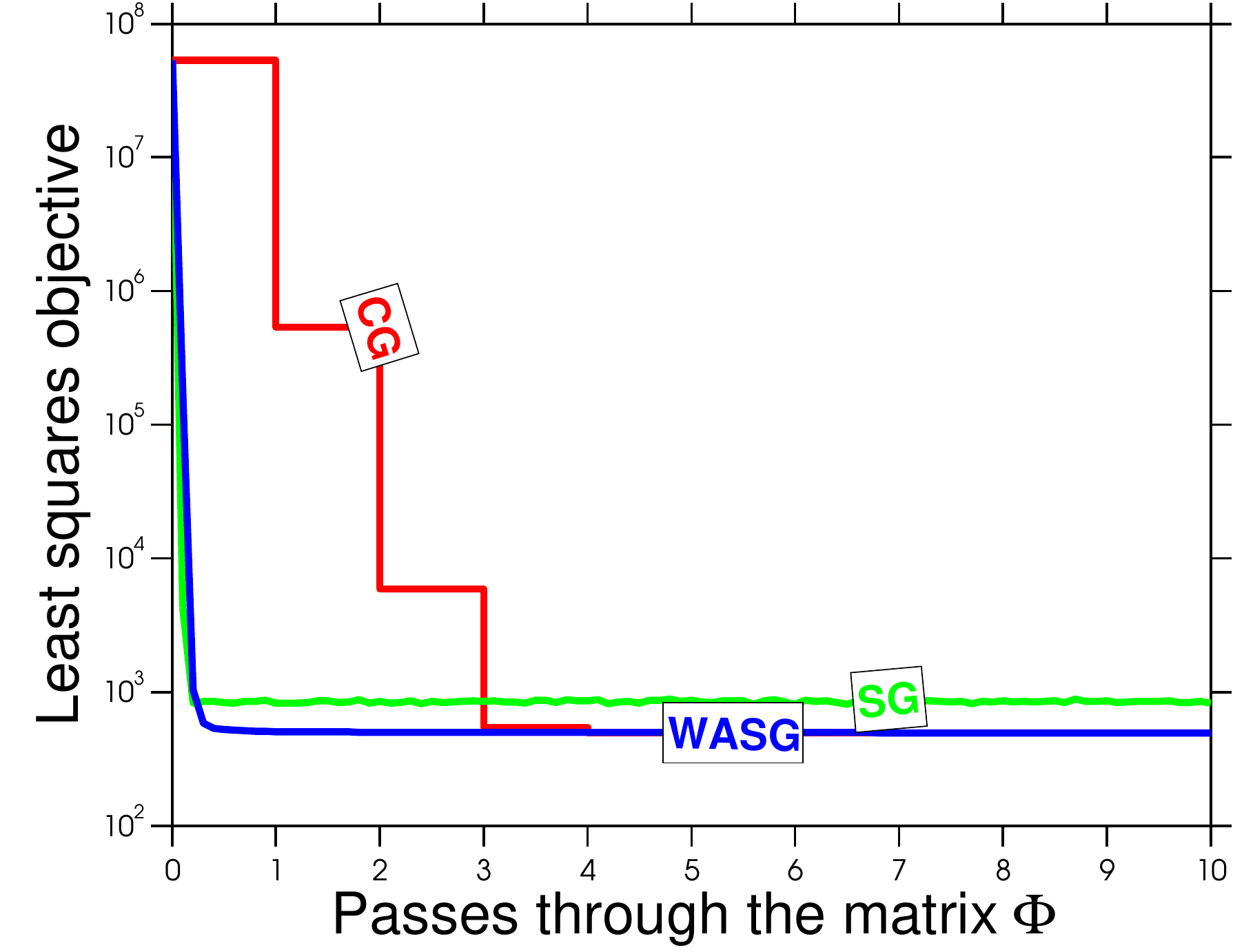}
    \includegraphics[width=0.35 \textwidth]{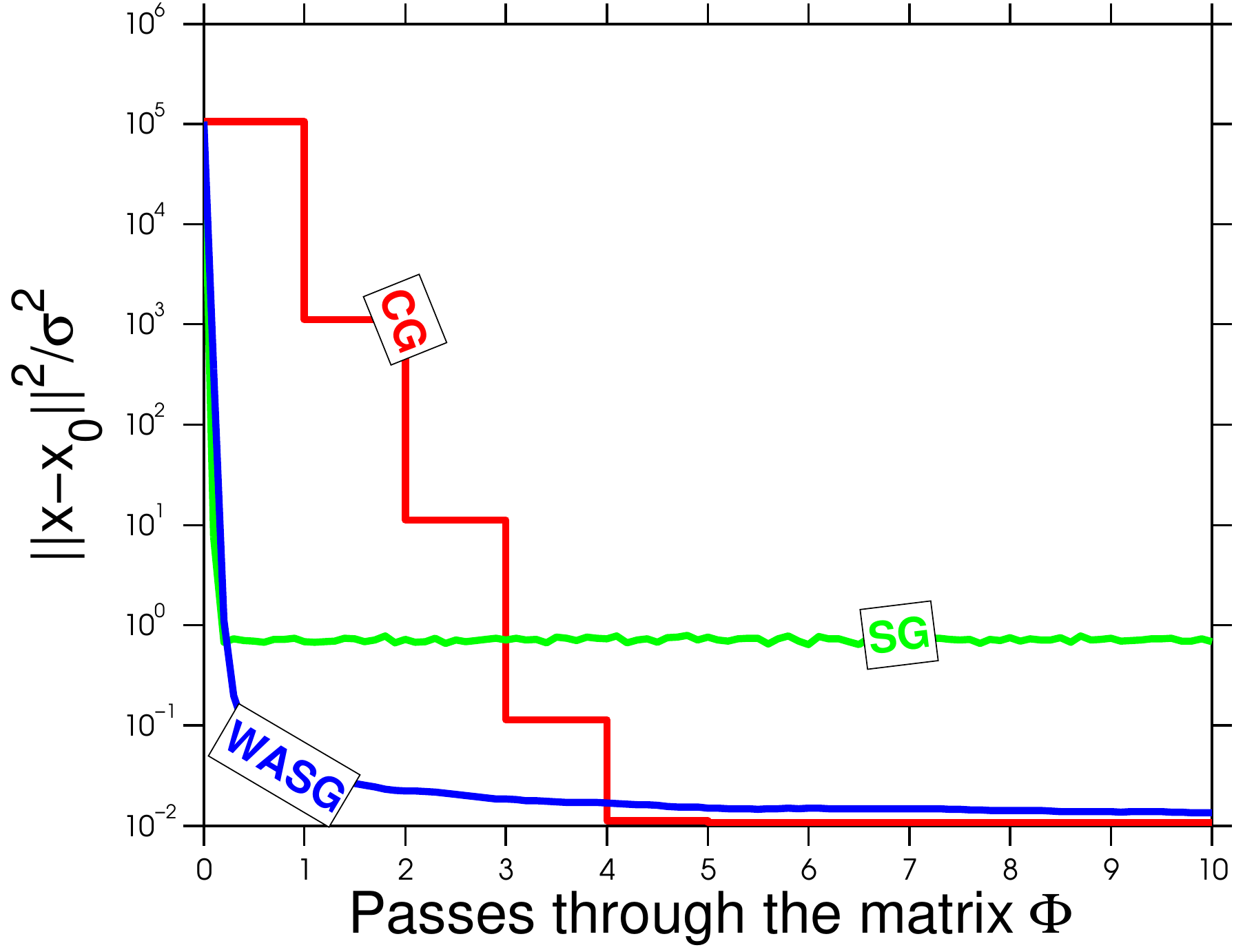}
  \caption{
  A numerical comparison of the CG method versus the stochastic gradient (SG) method and the weighted averaged SG (WASG) method for the LS problem \eqref{eq: LS}, showing the objective (a) and the normalized estimation error (b). The matrix $\Phi$ has standard normal entries with dimensions $n=10^5$ and $p=10^3$. The noise variance is $\sigma^2=10^{-2}$ whereas $\|x_0\|_2^2\approx p$. At a fractional access to the matrix $\Phi$, the stochastic methods obtain a good \emph{relative} accuracy on the signal estimate. Finally, the SG method has an optimization error due to our step-size choice; cf., Section II.B for an explanation.}\label{fig: intro sg} 
  %\vspace{-6mm}
\end{figure}

Figure \ref{fig: intro sg} shows that we can scale radically better than even the conjugate gradients (CG) method
for the LS formulation when $ n \gg p$  if we exploit stochastic approximation within first-order methods. We take the simplest optimization method, namely gradient descent with fixed step-size, and replace its gradient calculations with their cheap statistical estimates (cf., Section II for the recipe). 
The resulting \emph{stochastic} gradient algorithm already obtains a strong baseline performance with access to only a fraction of the rows of $\Phi$ while the conjugate gradient method requires many more full accesses.

 % spmag_smooth2
\section{First-Order Methods for Smooth and Non-Smooth Convex Optimization}
As the LASSO formulation highlights, non-smooth regularization can play an indispensable role in solution quality.  
By using the powerful \emph{proximal gradient} framework, we will see that many of these non-smooth problems can be solved nearly as efficiently as their smooth counterparts, a point not well-understood until the mid 2000s. To this end, this section describes first-order methods within this context, emphasizing specific algorithms with global convergence guarantees. In the sequel, we will assume that the readers have some familiarity with basic notions of convexity and complexity.  

%\vspace{-3mm}
\subsection{Smooth objectives}
We begin our exposition with an important special case of \eqref{eq: master}, where the objective $F$ only consists of a differentiable convex function $f$. The elementary first-order technique for this case is the  \emph{gradient} method, which uses only the local gradient $\nabla f(x)$ and iteratively performs the following update:
\begin{equation}\label{eq: grad method}
\xkk = \xk - \alpha_k\nabla f(\xk),
\end{equation}
where $k$ is the iteration count and $\alpha_k$ is an appropriate step-size that ensures convergence. 

For smooth minimization, we can certainly use several other \emph{faster} algorithms such as Newton-like methods.
By faster, we mean that these methods require fewer iterations than the gradient method to reach a target accuracy: 
 i.e., $F(\xk)- F^*\le \eps$. 
However, we do not focus on these useful methods since they either require additional information from the function $F$, more expensive computations, or do not generalize easily to constrained and non-smooth problems. 

Fortunately, the low per-iteration cost of the gradient method can more than make up for its drawbacks in iteration count. For instance, computing the gradient dominates the per-iteration cost of the method, which consists of  matrix-vector multiplications with $\Phi$ and its adjoint $\Phi^T$ when applied to the LS problem \eqref{eq: LS}. Hence, we can indeed perform \emph{many} gradient iterations for the cost of a single iteration of more complicated methods, potentially  taking a  shorter time to reach the same level of accuracy $\eps$.  

Surprisingly, by making simple assumptions about $f$, we can rigorously analyze how many iterations the gradient method will in fact need to reach an $\eps$-accurate solution. A common assumption that holds in many applications is that the gradient of $f$ is Lipschitz continuous, meaning that
\[
\forall x,y\in \mathbb{R}^p,\; \|\nabla f(x) - \nabla f(y)\|_2 \le L \|x-y\|_2, 
\]
for some constant $L$. 
When $f$ is twice-differentiable, a sufficient condition is the eigenvalues of its Hessian $\nabla^2 f(x)$ are bounded above by $L$. Hence, we can trivially estimate $L = \|\Phi\|_2^2$ for 
\eqref{eq: LS}.   

If we simply set the step-size $\alpha_k = 1/L$ or alternatively use a value that decreases $f$ the most, then the iterates of the gradient method for any convex $f$ with a Lipschitz-continuous gradient obey
\begin{equation}\label{eq: gradc}
f(\xk)-f^\star \le \frac{2L}{k+4}d_0^2,
\end{equation}
where $d_0 = \|\xzero - \xStar \|_2$ is the distance of the initial iterate $\xzero$ to an optimal solution $\xStar$
\cite[Cor.\ 2.1.2]{Nesterov2004}. Hence, the gradient method needs $\order(1/\eps)$-iterations for an $\eps$-accurate solution in the worst-case.

Unfortunately, this convergence rate  does not attain the known complexity lower-bound 
\[
f(\xk) - f^\star \ge \frac{3Ld_0^2}{32(k+1)^2},
\]
which holds for all functions $f$ with Lipschitz-continuous gradients. That is, in the worst case \emph{any} iterative method based only on function and gradient evaluations cannot hope for a better accuracy than $\Omega(1/k^2)$  at iteration $k$ for $k<p$~\cite{Nesterov2004}. 
Amazingly, a minor modification by Nesterov achieves this \emph{optimal} convergence by the simple step-size choice $\alpha_k = 1/L$ and an {extra}-momentum step with a parameter $\beta_k = \frac{k}{k+3}$ \cite{Nesterov2004}:
%%\vspace{-3mm}
\begin{algorithm}[H]
\begin{algorithmic}[1]
    \STATE $\xkk = \vk - \alpha_k\nabla f(\vk)$
    \STATE $\vkk = \xkk + \beta_k( \xkk - \xk)$
\end{algorithmic}
\caption{Nesterov's accelerated gradient method for unconstrained minimization ($\vzero=\xzero$) \cite{Nesterov2004}}
\label{alg:1}
\end{algorithm}
%\vspace{-3mm}
\noindent The \emph{accelerated} gradient method in Algorithm \ref{alg:1} achieves the best possible worst-case error rate, and hence, it is typically referred to as an \emph{optimal} first-order method. 

Many functions also feature additional structures useful for numerical optimization. Among them, \emph{strong convexity} deserves special attention since this structure provably offers key benefits such as the existence of a unique minimizer and improved optimization efficiency. A function $f$ is called strongly convex if the function $x\mapsto f(x) - \frac{\mu}{2}\|x\|_2^2$ is convex for some positive value $\mu$.  Perhaps not so obvious is the fact that even non-smooth functions can have strong convexity by this definition (i.e., $f(x) = \|x\|_1 + \mu/2\|x\|_2^2$).

Indeed, we can transform any convex problem into a strongly-convex problem by simply adding a squared $\ell_2$-regularization term. For instance, when we have $n<p$ in \eqref{eq: LS}, then the classic Tikhonov regularization results in a strongly convex objective with $\mu=\lambda$: 
\[
\widehat{x}_{\rm ridge} = \argmin_{x\in\Real^p} \left\{ F(x) :=  \frac{1}{2}\norm{y - \Phi x}_2^2 + \frac{\lambda}{2}\norm{x}_2^2\right\},
\]
The solution above is known as the ridge estimator  and offers statistical benefits \cite{wainwright2014structured}. When $f$ is twice-differentiable, a sufficient condition for strong convexity is that the eigenvalues of its Hessian $\nabla^2 f(x)$ are bounded below by $\mu$ for all $x$. For the LS problem \eqref{eq: LS}, strong convexity simply requires $\Phi$ to have independent columns. 

For strongly-convex problems with Lipschitz gradient, such as the ridge estimator, the gradient method geometrically converges to the unique minimizer when the step-size is chosen as $\alpha_k = 1/L$:
\begin{equation}\label{eq: sc}
\norm{\xk - \xStar}_2 \le \left(1 - \frac{\mu}{L}\right)^{k}\norm{\xzero - \xStar}_2.
\end{equation}
This convergence improves slightly when we instead use $\alpha_k=2/(\mu+L)$ \cite{{Nesterov2004}}.
Beside the obvious convergence rate difference, we highlight a subtlety between \eqref{eq: gradc} and \eqref{eq: sc}:  guarantees due to the Lipschitz assumption such as \eqref{eq: gradc} does not necessarily imply convergence in iterates $\xk$, while for stronlgy-convex functions we obtain guarantees on the convergence of both $f(\xk)$ and $\xk$.

It turns out that the accelerated-gradient method can also benefit from strong-convexity with an appropriate choice of the momentum term $\beta_k$. For example, if we set $\beta_k=(L-\mu)/(L+\mu)$, the accelerated gradient method obtains a near-optimal convergence rate given its assumptions~\cite[Thm.\ 2.2.3]{Nesterov2004}. In contrast, the gradient method automatically exploits strong convexity without any knowledge of $\mu$. 

\begin{table}[!t]\small
\begin{center}
\begin{tabular}{ccc}\toprule
 Algorithm & Convex & Strongly-Convex\\ 
 \midrule
\ [Proximal]-Gradient & $\order(Ld_0^2/\eps)$ & $\order\left(\frac{L}{\mu}\log(d_0^2/\eps)\right)$ \\
 Accelerated-[Proximal]-Gradient & $\order(\sqrt{Ld_0^2/\eps})$ & $\order\left(\sqrt{\frac{L}{\mu}}\log(d_0^2/\eps)\right)$\\
\bottomrule
\end{tabular}%\vspace{-6mm}
\end{center}
\caption{Total number of iterations to reach $\eps$-accurate solutions for first-order optimization methods. $L$ and $\mu$ denote the Lipschitz and strong convexity constants and $d_0= \|x^0-\xStar\|_2$. 
}
\label{table: taxonomy}
%\vspace{-6mm}
\end{table}

\begin{figure}[th]
\begin{center}
    \includegraphics[height=2in]{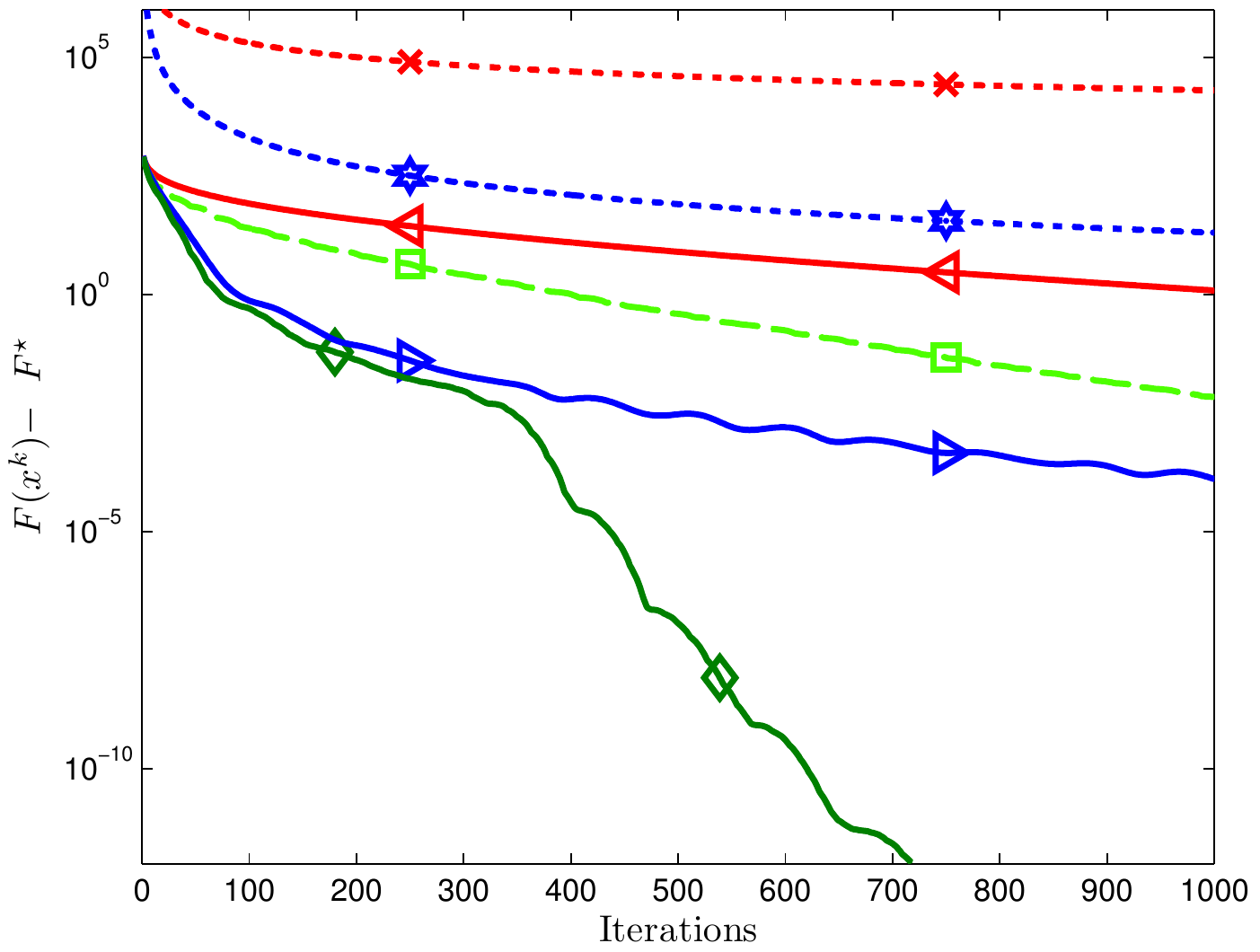}
    \includegraphics[height=2in]{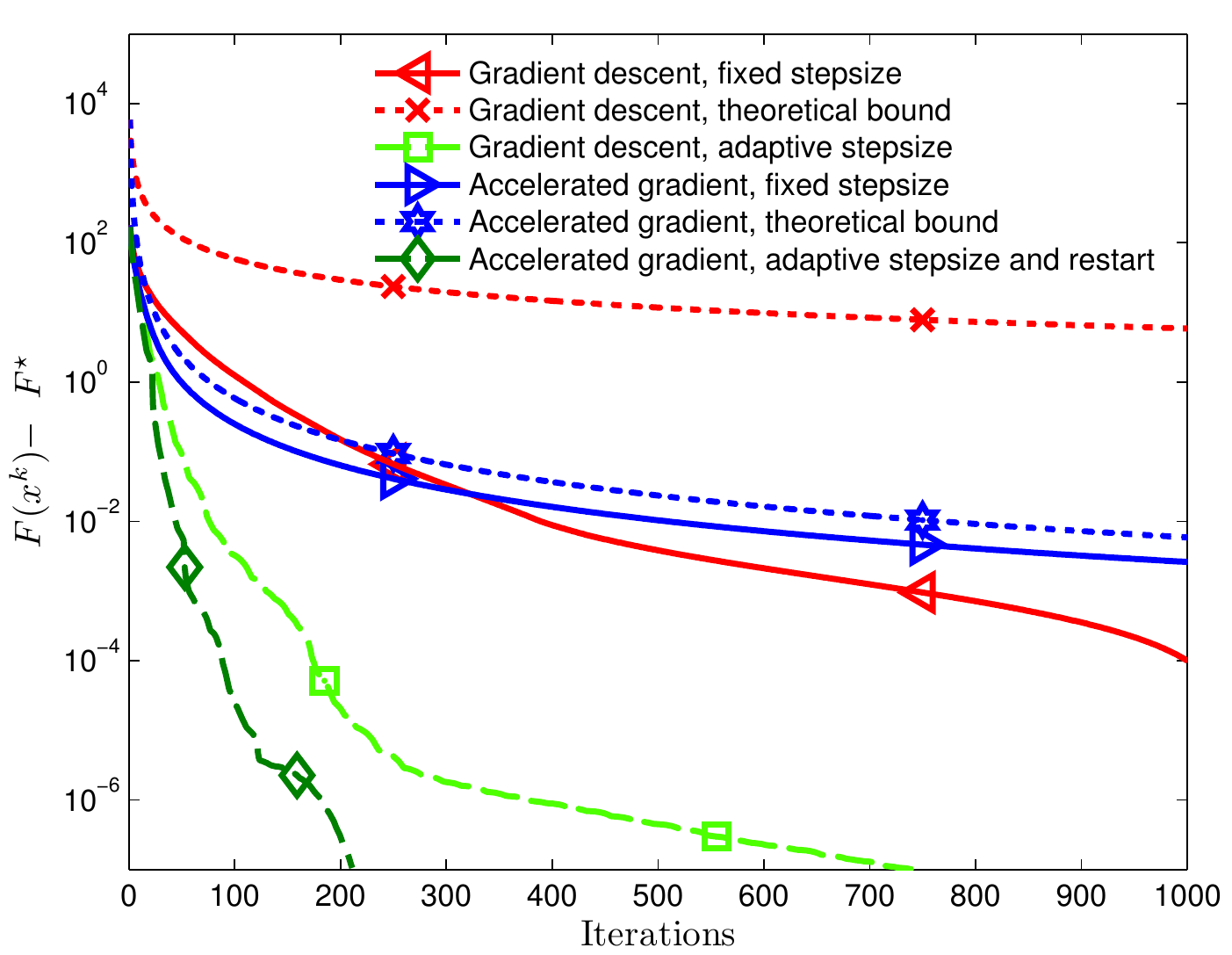}
\end{center}
%\vspace{-1em}
\caption{
The performance of first-order methods can improve significantly with practical enhancements. We demonstrate how the objective 
$F(\xk)$ 
progresses as a function of iterations $k$ 
for solving (Left) the LS formulation \eqref{eq: LS}, and (Right) the LASSO formulation \eqref{eq: lasso}, both with 
$p=5000$ and $n=2500$,    
for four methods: (proximal)-gradient descent with fixed step-size $\alpha=1/L$ and adaptive step-size, accelerated (proximal)-gradient descent with fixed step-size, and accelerated (proximal)-gradient descent with the adaptive step-size and restart scheme in TFOCS~\cite{TFOCS}. For the LS formulation, the basic methods behave qualitatively the same as their theoretical upper-bounds predict but dramatically improve with the enhancements. 
In the LASSO formulation, gradient descent automatically benefits from sparsity of the solution and actually outperforms the basic accelerated method in high-accuracy regime, but adding the adaptive restart enhancement allows the accelerated method to also benefit from sparsity. 
}
\label{fig:enhancements}
%\vspace{-6mm}
\end{figure}

Table~\ref{table: taxonomy} summarizes the number of iterations to reach an accuracy of $\epsilon$ for the different configurations discussed in this section. Note however that there are numerous practical enhancements, such as step-size selection rules for $\alpha_k$ and adaptive restart of the momentum parameter $\beta_k$ \cite{adaptiveRestart} that add only a small computational cost and do not rely on knowledge of the Lipschitz constant $L$ or the strong-convexity parameter $\mu$. 
While such tricks-of-the-trade do not rigorously improve the worst-case convergence rates, they often lead to superior empirical convergence (cf., Figure \ref{fig:enhancements}) and similarly apply to their important \emph{proximal} counterparts for solving \eqref{eq: master} that we discuss next.

Finally, the fast gradient algorithms described here also apply to non-smooth minimization problems using Nesterov's smoothing technique \cite{Nesterov2004}. In addition,  rather than assuming Lipschitz-continuity of the gradient of the objective function, recent work has considered efficient gradient methods for smooth self-concordant functions, which naturally emerge in Poisson imaging, graph learning, and quantum tomography problems~\cite{tran2013composite}.

 % spmag_proximal
%\vspace{-3mm}
\subsection{Composite objectives}
We now consider the canonical composite problem \eqref{eq: master}, where the objective $F$ consists of a differentiable convex function $f$ and a non-smooth convex function $g$ as in~\eqref{eq: lasso}. 

In general, the non-differentiability of $g$ seems to substantially reduce the efficiency of first-order methods. This was indeed the conventional wisdom since generic non-smooth optimization methods, such as subgradient and bundle methods, require $\mathcal{O}(1/\eps^2)$ iterations to reach $\eps$-accurate solutions~\cite{Nesterov2004}. While strong convexity helps to improve this rate to $\mathcal{O}(1/\eps)$, the resulting rates are slower than first-order methods for smooth objective.

Fortunately, composite objectives are far from generic non-smooth convex optimization problems. The \emph{proximal-gradient} methods  specifically take advantage of the composite  structure in order to retain the same convergence rates of the gradient method for the smooth problem classes in Table \ref{table: taxonomy}~\cite{SchmidtNIPS11}. 
It becomes apparent that these algorithms are in fact natural extensions of the gradient method when we view the gradient method's iterations \eqref{eq: grad method} as an optimization problem: 
\begin{equation}\label{eq: grad iter}
\xkk = \argmin_{y\in \R^p} \left\{ f(\xk) + \nabla f(\xk)^T(y-\xk) + \frac{1}{2\alpha_k}\norm{y-\xk}^2\right\},
\end{equation}
which is based on a simple local quadratic approximation of $f$. Note that when $\alpha_k \leq 1/L$, the objective function above is a quadratic upper bound on $f$. 
Proximal-gradient methods use the same approximation of $f$, but simply include the non-smooth term $g$ in an explicit fashion:
\begin{equation}\label{eq: prox eq}
\xkk = \argmin_{y\in \R^p} \left\{f(\xk) + \nabla f(\xk)^T(y-\xk) + \frac{1}{2\alpha_k}\norm{y-\xk}^2 + g(y)\right\}.
\end{equation}
For $\alpha_k \leq 1/L$, the objective is an upper bound on $F$ in \eqref{eq: master}.

The optimization problem \eqref{eq: prox eq} is the update rule of the proximal-gradient method:
\[
\xkk = \prox_{\alpha_kg}(\xk - \alpha_k\nabla f(\xk)),
\]
where the \emph{proximal map} or \emph{proximal operator} is defined as
\begin{equation}\label{eq: prox}
  \prox_{g}(y) \defeq \argmin_{x} \left\{ g(x)+ \frac{1}{2}\|x-y \|_2^2 \right\}.
\end{equation}
The \emph{accelerated} proximal-gradient method is defined analogously: 
%\vspace{-3mm}
\begin{algorithm}[H]
\begin{algorithmic}[1]
    \STATE $\xkk = \prox_{\alpha_k g}\left( \vk - \alpha_k\nabla f(\vk) \right)$
    \STATE $\vkk = \xkk + \beta_k( \xkk - \xk)$
\end{algorithmic}
\caption{Accelerated proximal gradient method to solve \eqref{eq: master}~\cite{BecTeb:09,Nesterov2004}. Set $\vzero=\xzero$.
}
\label{alg:2}
\end{algorithm}
%\vspace{-3mm}
  
An interesting special case of the proximal-gradient algorithm arises if we consider the indicator function on a convex set $\mathcal{C}$, which is an elegant way of incorporating constraints into \eqref{eq: master} 
\[\small
g(x) = \begin{cases}
0 & x \in \mathcal{C}\\
\infty & x \notin \mathcal{C}
\end{cases}.
\]
Then, the proximal-gradient method yields the classic projected-gradient method for constrained optimization.
  
These method's fast convergence rates can also be preserved under approximate proximal maps~\cite{SchmidtNIPS11}.
Proximal operators offer a flexible computational framework to incorporate a rich set of signal priors in optimization. For instance, we can often represent a signal $x_0$ as a linear combination of atoms  $ a \in \mathcal{A}$ from some \emph{atomic set} $\mathcal{A}\subseteq \R^p$ as $x_0 = \sum_{a \in \mathcal{A}} c_{a} a,$ where $c_a$ are the representation coefficients. Examples of atomic sets include structured sparse vectors, sign-vectors, low-rank matrices, and many more. The geometry of these sets can facilitate perfect recovery even from underdetermined cases of the linear observations \eqref{eq: linear model} with sharp sample complexity characterizations \cite{chandrasekaran2012convex}.

To promote the structure of the set  $\mathcal{A}$ in convex optimization, we can readily exploit its gauge function: $ g_{\mathcal{A}}( x) \defeq \inf \big \{\rho > 0 \mid x \in \rho \cdot \overline{\rm conv}(\mathcal{A}) \big \}$, where $\overline{\rm conv}(\mathcal{A})$ is the convex hull of the set $\mathcal{A}$. The corresponding proximal operator of the gauge function has the following form 
\begin{equation}\label{eq: atomic prox}
  \mathrm{prox}_{\gamma g_\mathcal{A}} (u) = u  - \argmin_{ v \in\R^d}\left \{\| u -  v \|^2_2:  \langle  a,   v \rangle \le \gamma, \forall   a \in \mathcal{A} \right\},
\end{equation}
which involves a quadratic program in general but can be explicitly calculated in many cases \cite{chandrasekaran2012convex}. Intriguingly, \eqref{eq: atomic prox} also has mathematical connections to \emph{discrete} submodular minimization. 

By and large, whenever the computation of the proximal map is efficient, so are the proximal-gradient algorithms. For instance, when $g(x) = \lambda\|x\|_1$ as in the LASSO formulation \eqref{eq: lasso}, the proximal operator is the efficient soft thresholding operator. Against intuition, a set with an infinite number of atoms can admit an efficient proximal map, such as the set of rank-1 matrices with unit Frobenius norm whose proximal map is given by singular value thresholding. On the other hand, a set with a finite number of atoms need not, such as rank-1 matrices with $\pm 1$ entries whose proximal operator is intractable. Numerous other examples exist \cite{CombettesPesquetChapter,chandrasekaran2012convex}.

\begin{figure}[ht]
\begin{center}
\includegraphics[width=3in]{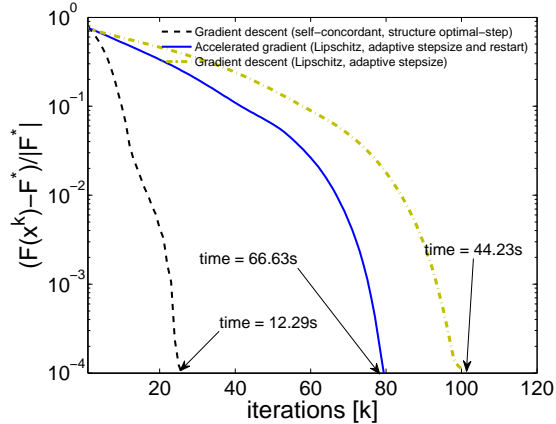}
\end{center}
%\vspace{-1em}
\caption{ Choosing the correct smoothness structure on $f$ in composite minimization is key to the numerical efficiency of the first-order methods. The convergence plot here simply demonstrates this with a composite objective, called the heteroskedastic LASSO (hLASSO) from \cite{tran2013composite} where $n=1.5\times 10^4$ and $p=5\times 10^4$. hLASSO features  a self-concordant but not Lipschitz gradient smooth part $f$, and obtains sparse solutions while simultaneously estimating the unknown noise variance $\sigma^2$ in the linear model \eqref{eq: linear model}. The simplest first-order method, when matched to correct self-concordant smoothness structure, can calculate its step-sizes optimally, and hence, significantly outperforms the standard first-order methods  based on the Lipschitz gradient assumption even with the enhancements we discussed. Surprisingly, the accelerated gradient method takes longer than the gradient method to reach the same accuracy since it relies heavily on the Lispchitz gradient assumption in its momentum steps, which lead to costlier step-size adaptation.}
\label{fig: selfcon}
%\vspace{-6mm}
\end{figure}

When $g$ represents the indicator function of a compact set, the Frank-Wolfe method solves~\eqref{eq: prox eq} without the quadratic term and can achieve an $O(1/\epsilon)$ convergence rate in the convex case~\cite{jaggi13}. This linear sub-problem may be easier to solve, and the method gives sparse iterations since each method only modifies a single element of the atomic set.

Finally, proximal-gradient methods that optimally exploit the self-concordance properties of $f$  are also explored in \cite{tran2013composite} (c.f., Figure \ref{fig: selfcon} for an example). Interestingly, many self-concordant functions themselves have tractable proximity operators as well---a fact that proves useful next.

 % spmag2_admm
%%\vspace{-3mm}
\subsection{Proximal objectives}

For many applications, the first-order methods we have covered so far are not directly applicable. As a result,  we will find it useful here to view the composite form \eqref{eq: master} in the following guise 
\begin{equation}\label{eq:Douglas-Rachford}
\min_{x,z\in \mathbb{R}^p}\left\{F(x,z) \defeq h(x) + g(z): \Phi z = x\right\},
\end{equation}
and only posit that the proximity operators of $h$ and $g$ are both efficient. 

This seemingly innocuous reformulation can simultaneously enhance our modeling and computational capabilities. First, \eqref{eq:Douglas-Rachford} can address non-smooth and non-Lipschitz objective functions that commonly occur in many applications  \cite{mccoy2013convexity,tran2013composite}, such as robust principal component analysis (RPCA), graph learning, and Poisson imaging, \emph{in addition} to the composite objectives we have covered so far. Second, we can apply a simple algorithm, called the alternating direction method of multipliers (ADMM) for its solutions, which leverages powerful augmented Lagrangian and dual decomposition techniques\cite{CombettesPesquetChapter,boyd2011distributed}:
%\vspace{-3mm}
\begin{algorithm}[H]
\begin{algorithmic}[1]
    \STATE $\xkk = \argmin_x\; \gamma h(x) + \frac{1}{2}\|x - \Phi\zk+\uk\|_2^2 = \prox_{\gamma h}( \Phi\zk - \uk)$  %
    \STATE $\zkk = \argmin_z\; \gamma g(z) + \frac{1}{2}\|\xkk - \Phi z + \uk\|_2^2 $
    \STATE $\ukk = \uk + \xkk - \Phi\zkk$
\end{algorithmic}
\caption{ADMM to solve \eqref{eq:Douglas-Rachford}; $\gamma>0, \zzero=\uzero=0$}
\label{alg:3}
\end{algorithm}
%\vspace{-3mm}

Algorithm \eqref{alg:3} is well-suited for distributed optimization and turns out to be equivalent or closely related to many other algorithms, such as Douglas-Rachford splitting and Bregman iterative algorithms \cite{boyd2011distributed}. ADMM requires a penalty parameter $\gamma$ as input and produces a sequence of iterates that approach feasibility and produce the optimal objective value in the limit. An overview of ADMM, its convergence, enhancements, parameter selection, and stopping criteria can be found here  \cite{boyd2011distributed}.

We highlight two {caveats} for ADMM. First, we have to numerically solve step 2 in  Algorithm \ref{alg:3} in general except when $\Phi^T\Phi$ is efficiently diagonalizable.
Fortunately, many notable applications support these features, such as matrix completion where $\Phi$ models sub-sampled matrix entries, image deblurring where $\Phi$ is a convolution operator, and total variation regularization where $\Phi$ is a differential operator with periodic boundary conditions. 
Secondly, the \emph{na\"ive} extension of ADMM to problems with more than two objective terms no longer has convergence guarantees.

%\vspace{-3mm}
\begin{algorithm}[H]
\begin{algorithmic}[1]
    \STATE $\xkk = \prox_{\gamma h}( \Phi\zk - \uk)$  %
    \STATE $\zkk = \prox_{\gamma\tau g}( \zk + \tau\Phi^T(\xkk-\Phi\zk+\uk))$
    \STATE $\ukk = \uk + \xkk - \Phi\zkk$
\end{algorithmic}
\caption{Primal-Dual Hybrid Gradient algorithm to solve \eqref{eq:Douglas-Rachford}; $\gamma>0$ and $\tau \le 1/\|\Phi\|^2$}
\label{alg:lADMM}
\end{algorithm}
%\vspace{-3mm}
 
Several solutions address the two drawbacks above. For the former, we can update $\zkk$ inexactly by using a single step of the proximal gradient method, 
which leads to the method shown in Algorithm~\ref{alg:lADMM} 
which was motivated in \cite{preconditionedADMM} as a preconditioned variant of ADMM and then analyzed in \cite{ChambollePock10} in a more general framework.
Interestingly, when $h(x)$ has a difficult proximal operator in Algorithm \ref{alg:3} but also has a Lipschitz gradient, we can replace $h(x)$ in Step 1 with its quadratic surrogate as in \eqref{eq: grad iter} to obtain the \emph{linearized} ADMM \cite{Condat2011}. 
Surprisingly, these inexact update-steps can be as fast to converge as the full ADMM in certain applications \cite{preconditionedADMM}. 
We refer the readers to \cite{Condat2011,boyd2011distributed,ChambollePock10} for the parameter selection of these variants as well as their convergence and generalizations.

For the issue regarding objectives with more than two terms, we can use dual decomposition techniques to treat the multiple terms in the objective of \eqref{eq:Douglas-Rachford} as separate problems and simultaneously solve them \emph{in parallel}. We defer  this to Section III and Algorithm~\ref{alg:4}.

 % spmag_scaling_up
\section{Big Data scaling via randomization}
In theory, first-order methods are well-positioned to address very large-scale problems. \emph{In practice}, however, the exact numerical computations demanded by their iterations can make even these simple methods infeasible as the problem dimensions grows. Fortunately, it turns out  that first-order methods are quite robust to using approximations of their optimization primitives, such as gradient and proximal calculations\cite{SchmidtNIPS11}.
This section describes emerging randomized approximations that increase the reach of first-order methods to extraordinary scales.  

To deliver an example of the key ideas, we will focus only on the smooth and strongly convex $F$ as objectives and point out extensions when possible. Many notable Big Data problems indeed satisfy this assumption. For instance, Google's PageRank problem 
measures the importance of nodes in a given graph via its incidence matrix $M \in \mathbb{R}^{p\times p}$ and $p$ is on the order of \emph{tens of billions}. Assuming that more important nodes have more connections, the problem in fact aims to find the top singular vector of the {stochastic} matrix $\Phi = M {\rm diag}(M^T\mathbf{1}_p)^{-1}$, where $\mathbf{1}_p\in \mathbb{R}^p$ is the vector of all $1$'s. 

The  PageRank algorithm simply solves this basic linear algebra problem (i.e., find $x^*\ge 0$ such that $\Phi x^* = x^*$ and $\mathbf{1}_p^Tx^*=1$) with the \emph{power method}. However, we can well-approximate this goal using a least squares problem when we relax the constraints with a penalty parameter $\gamma>0$ \cite{nesterov2012efficiency}:
\begin{equation}\label{eq: google}
  \min_{x \in \mathbb{R}^p}\left \{ F(x)\defeq \frac{1}{2} \|x- \Phi x\|_2^2 +\frac{\gamma}{2}\left( \mathbf{1}_p^Tx-1\right)^2 \right\},
\end{equation}
Note that we can work with a constrained version of this problem that includes positivity, but since the PageRank formulation itself  is not exact model of reality, the simpler problem can be preferable for obvious computational reasons.  
Clearly, we would like to minimize the number of operations involving the matrix $\Phi$ in any solution method. 

%\vspace{-3mm}
\subsection{Coordinate descent methods}
Calculating the full gradient for the PageRank problem formulation requires a matrix-vector operation at each iteration. A cheaper vector-only operation would be to pick a coordinate $i$ of $x$ and only modify the corresponding variable $x_i$ to improve the objective function.
This idea captures the essence of \emph{coordinate descent methods}, which have a long history in optimization~\cite{luo1992convergence} and are related to classic methods like 
the Gauss-Seidel cyclic reduction strategy for solving linear systems. The general form of coordinate descent methods is illustrated in Algorithm~\ref{alg:coordinateDescent}, where $e_i$ is the $i^\text{th}$ canonical coordinate vector and $\nabla_{i}F(\cdot)$ is the $i^\text{th}$ coordinate of the gradient. 

%\vspace{-3mm}
\begin{algorithm}[H]
\begin{algorithmic}[1]
    \STATE Choose an index $i_k\in\{1,2,\ldots,p\}$ \emph{(see the main text for possible selection schemes)}
    \STATE $\xkk = \xk - \alpha\nabla_{i_k}F(\xk)e_{i_k}$
\end{algorithmic}
\caption{Coordinate descent to minimize $F$ over $\R^p$}
\label{alg:coordinateDescent}
\end{algorithm}
%\vspace{-3mm}

The key design consideration across all coordinate descent methods is the choice of the coordinate $i$ at each iteration. A simple strategy amenable to analysis is to greedily pick the coordinate with the largest \emph{directional derivative} $\nabla_i F$. This selection with $\alpha=1/L_\text{max}$ or optimizing the variable exactly leads to a convergence rate of
\begin{equation}
\label{eq:CD-rate}
 F(\xk) - F(x^*) \leq \left(1 - \frac{\mu}{pL_{\text{max}}}\right)^k(F(\xzero) - F(\xStar)),
\end{equation} 
where $L_{\text{max}} \defeq \max_i L_i$ is the maximum across the Lipschitz constants of $\nabla_{i}F(x)$~\cite{nesterov2012efficiency}. This configuration indeed seeks the best reduction in the objective per iteration  we can hope for under this setting. 

The example above underlines the fundamental difficulty in coordinate descent methods. Finding the best coordinate to update, the maximum of the gradient element's magnitudes, can require a computational effort as high as the gradient calculation itself. However, the incurred cost is not justified since the method's convergence is provably slower than the gradient method due to the basic relationship $L_i \leq L \leq pL_i$. An alternative proposal 
is to cycle through all coordinates sequentially. This is the cheapest coordinate selection strategy we can hope for 
but it results in a substantially slower convergence rate. 

Surprisingly, randomization of the coordinate choice can achieve the best of both worlds.
Suppose we choose the coordinate $i$ uniformly at random among the set $\{1,2,\dots,p\}$. This selection can be done with a cost independent of $p$, but surprisingly nevertheless achieves the same convergence rate~\eqref{eq:CD-rate} in expectation~\cite{nesterov2012efficiency}. The randomized algorithm's variance around its expected performance is well-controlled.

We also highlight two salient features of coordinate descent methods. First, they are perhaps most useful for objectives of the form $F(Ax)$ with $A \in \mathbb{R}^{n\times p}$, where evaluating the (not necessarily smooth) $F$  costs $\order(n)$. By tracking the product $A\xk$ with incremental updates, we can then perform coordinate descent updates in linear time. 
Second, if we importance sample the coordinates proportional to their Lipschitz constants $L_i$, then the convergence rate of the randomized method improves to
\begin{equation}
\label{eq:CD-rate2}
 F(\xk) - F(\xStar) \leq \left(1 - \frac{\mu}{pL_{\text{mean}}}\right)^k(F(\xzero) - F(\xStar)),
\end{equation}
where $L_{\text{mean}}$ is the mean across the $L_i$. Hence, this non-uniform random sampling strategy improves the speed by only adding an  $\order{(\log(p))}$ importance sampling cost to the algorithm  \cite{nesterov2012efficiency}. 

Finally, accelerated and composite versions of coordinate descent methods have also recently been explored, although accelerated methods often do not preserve the cheap iteration cost of the non-accelerated versions \cite{nesterov2013first}: cf., \cite{nesterov2012efficiency} for a numerical example of these methods on the PageRank problem \eqref{eq: google}. 

%\vspace{-3mm}
\subsection{Stochastic gradient methods}
In contrast to randomized coordinate descent methods, which update a single coordinate at a time with its exact gradient, stochastic gradient methods update all coordinates simultaneously but use approximate gradients. They are best suited for minimizing decomposable objective functions $F$
\begin{equation}
\label{eq:SG}
\min_{x\in \mathbb{R}^p}\left\{ F(x) \defeq  \frac{1}{n}\sum_{j=1}^n F_j(x)\right \},
\end{equation}
where each $F_j$ measures the data misfit for a single data point. This includes models as simple as least squares and also more elaborate models like conditional random fields. 

%\vspace{-3mm}
\begin{algorithm}[H]
\begin{algorithmic}[1]
    \STATE Choose an index $j_k\in \{1,2,\ldots,n\}$ \emph{uniformly at random}
\STATE $\xkk = \xk - \alpha_k \nabla F_{j_k}(\xk)$ 
\end{algorithmic}
\caption{Stochastic gradient descent to minimize $F$ over $\R^p$}
\label{alg:stochasticGradient}
\end{algorithm}
%\vspace{-3mm}

Stochastic gradient iterations heavily rely on the decomposability of~\eqref{eq:SG} as shown in Algorithm \ref{alg:stochasticGradient}.
Similar to the coordinate descent methods, the crucial design problem in the stochastic gradient methods is the selection of the data points $j$ at each iteration. Analogously, we obtain better convergence rates by choosing the $j$ uniformly at random rather than cycling through the data~\cite{nedic2001convergence}. In contrast, per iteration cost of the algorithm now depends only on $p$ but not $n$.

Interestingly, a random data point selection results in an unbiased gradient estimate when we view \eqref{eq:SG} as the empirical observation of a  \emph{expected risk function} that governs the optimization problem
\begin{equation}\label{eq: true cost}
\min_{x\in \mathbb{R}^p}\left\{ F(x) \defeq  \mathbb{E}_\xi[ F_\xi(x)]\right\},
\end{equation}
where the expectation is taken over the sampling distribution for the indices $\xi$. Indeed, if we can sample from the true underlying distribution, then stochastic gradient methods directly optimize the expected risk minimization problem and  result in provable generalization capabilities in machine learning applications \cite{bottou2007tradeoffs}. 
The unbiased gradient estimation idea enables the stochastic gradient descent method to handle convex objectives beyond the decomposable form \eqref{eq:SG} \cite{nesterov2013first}.

The general SG method has classically used a 
 decreasing sequence of step-sizes $\{\alpha_k\}$. However, this unfortunately leads to the same slow $\mathcal{O}(1/\sqrt{\epsilon})$ and $\mathcal{O}(1/\epsilon)$ convergence rates of the sub-gradient method.
But interestingly, if we still use a constant step-size at each iteration for the stochastic gradient method the algorithm is known to quickly reduce the initial error, even if it has a non-vanishing optimization error~\cite{nedic2001convergence}. We have observed this for the stochastic gradient descent example in Figure \ref{fig: intro sg}. 
Indeed,
while stochastic gradient descent methods have historically been notoriously hard to tune, recent results show that using large step sizes and weighted averaging of the iterates (cf., Figure \ref{fig: intro sg}) allows us to achieve optimal convergence rates while being robust to the setting of the step size and the modeling assumptions~\cite{nedic2001convergence,bach2013nonStrongly}.
For example, recent work has shown~\cite{bach2013nonStrongly} that an averaged stochastic gradient iteration with a constant step-size achieves an $\order(1/\eps)$ convergence rate even without strong convexity under joint self-concordance-like and Lipschitz gradient assumptions. Another interesting recent development has been stochastic algorithms that achieve linear convergence rates for strongly-convex problems of the form~\eqref{eq:SG} in the special case where the data size is finite~\cite{schmidt2013finite}.

 %
 
 %\vspace{-3mm} 
\subsection{Randomized linear algebra}
\newcommand{\XRLA}{M}
For Big Data problems, basic linear algebra operations, such as matrix decompositions (e.g., eigenvalue, singular value, and Cholesky) and matrix-matrix  multiplications can be major computational bottlenecks due to their superlinear dependence on dimensions.
However, when the relevant matrix objects have low-rank representations (i.e., $\XRLA = LR^T$ with $L\in \R^{p\times r}$ and $L\in \R^{p\times r}$ where $r \ll p$), the efficiency of these methods uniformly improves. For instance, the corresponding singular value decomposition (SVD) of $\XRLA$ would only cost $\order(pr^2 + r^3)$ flops. %

The idea behind randomized linear algebra methods is either to approximate $\XRLA \approx Q(Q^T\XRLA)$ with $Q\in\R^{p\times r}$, or to construct a low-rank representation by column or row subset selection in order to speed up computation.  And indeed, doing this in a \emph{randomized} fashion gives us control over the distribution of the errors~\cite{halko2011finding,MahoneyMonograph}.  This idea generalizes to matrices of any dimensions and  has the added benefit of exploiting mature computational routines in  nearly all programming languages. Hence, they immediately lend themselves well to modern distributed architectures. 

We describe three impacts of randomizing linear algebra routines in optimization here.  First, we can accelerate computation of the proximity operators of functions that depend on spectral values of a matrix. 
For instance, the proximity operator of the nuclear norm, used in matrix completion and RPCA problems, requires a partial SVD. This is traditionally done with the Lanczos algorithm which does not parallelize easily due to synchronization and re-orthogonalization issues. However, with the randomized approach, the expected error in the computation is bounded and can be used to maintain rigorous guarantees for the convergence of the whole algorithm~\cite{becker2013randomized}.

Secondly, the idea also works in obtaining unbiased gradient estimates for matrix objects, when randomization is chosen appropriately, and hence applies to virtually all stochastic gradient algorithms. 
Finally, the randomized approach can be used to \emph{sketch} objective functions, i.e., to approximate them in order to obtain much cheaper iterations with exact first-order methods while retaining accuracy guarantees for the true objective~\cite{MahoneyMonograph}.

%\vspace{-3mm}
\begin{algorithm}[H]
\begin{algorithmic}[1]
	\REQUIRE $\XRLA \in\R^{p \times p}$, integer $r$
    \STATE Draw $\Omega \in \R^{p \times r}$ iid $\mathcal{N}(0,1)$
    \STATE $W =\XRLA \Omega$ \COMMENT{Matrix multiply, cost is $\order(p^2r)$}
		\STATE $QR=W$ \COMMENT{QR algorithm, e.g., Gram-Schmidt, cost is $\order(pr^2)$}
		\STATE $U=\XRLA^TQ$ \COMMENT{Matrix multiply, cost is  $\order(p^2r)$}
	\RETURN $\widehat{\XRLA}_{(r)} = QU^T$ \COMMENT{Rank $r$}
\end{algorithmic}
\caption{Randomized low-rank approximation
}
\label{alg:svd}
\end{algorithm}
%\vspace{-3mm}

Algorithm~\ref{alg:svd} is an example of a randomized low-rank approximation, which is simply a single step of the classical QR iteration, using a random initial value. Surprisingly, the error in approximating $\XRLA$ is nearly as good as the \emph{best} rank-$r$ approximation, where $\ell=r+\rho$ and $\rho$ is small. Specifically, for $r\ge 2, \rho \ge 2$ and $\ell \le p$, \cite{halko2011finding} provides the bound
\[
\E \|\widehat{\XRLA}_{(\ell)}-\XRLA\|_F \le \sqrt{1+\frac{r}{\rho-1}}\|\XRLA- \XRLA_{(r)}\|_F
\]
where the expectation is taken with respect to the randomization, $\XRLA_{(r)}$ is the best rank-$r$ approximation of $\XRLA$, which only keeps the first $r$ terms in the SVD and sets the rest to zero; furthermore, \cite{halko2011finding} shows a deviation bound showing that the error concentrates tightly around the expectation. Thus, the approximation can be very accurate if the spectrum of the matrix decays to zero rapidly. For additional randomized linear algebra schemes and their corresponding guarantees, including using a power iteration to improve on this bound, we refer the readers to \cite{halko2011finding}. 

\begin{figure}[ht]
\centering
\includegraphics[width=2.6in]{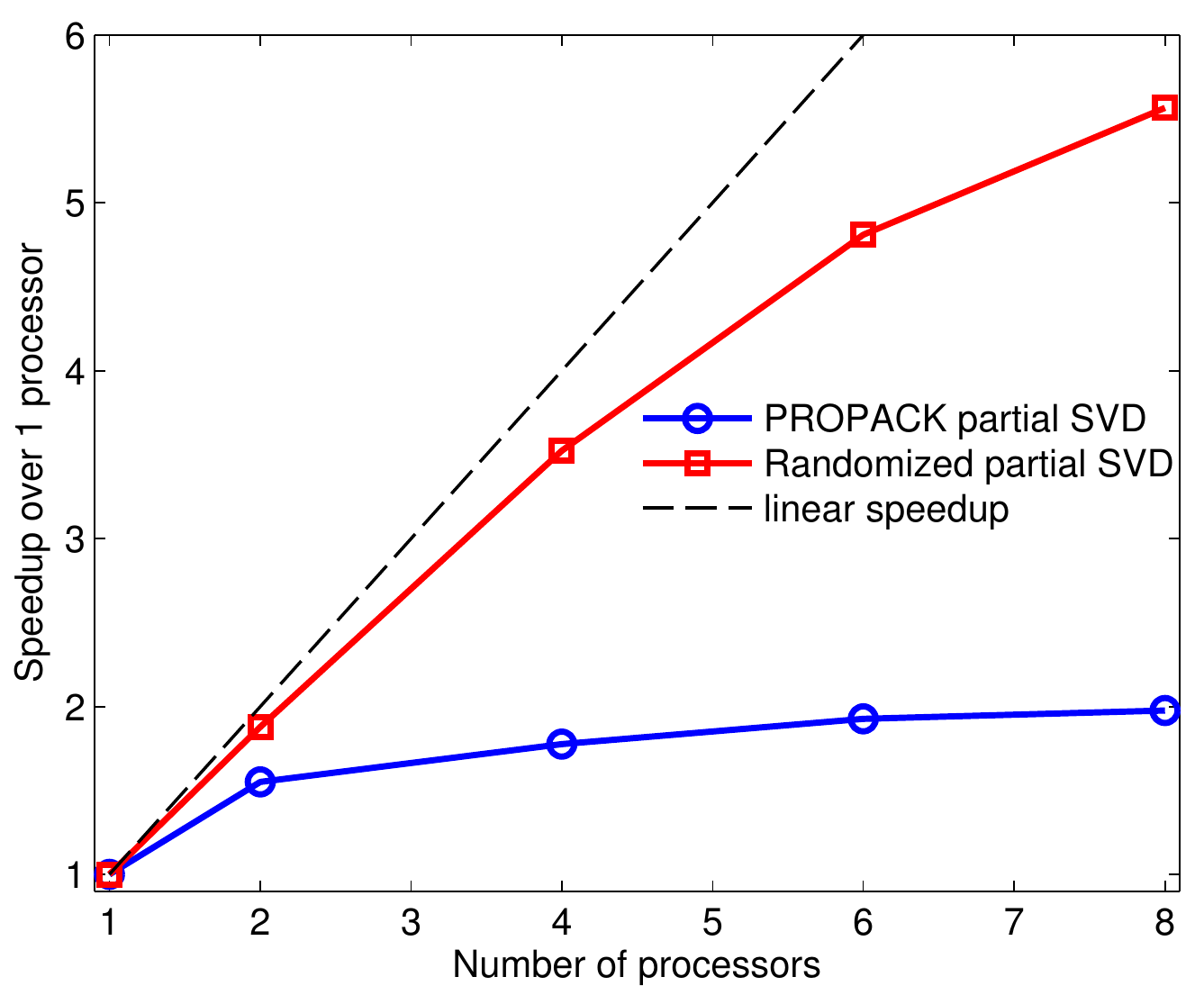}
%\vspace{-1em}
\caption{
Computing the top 5 singular vectors of a $10^9$ entry matrix using varying number of computer cores. 
The matrix is a dense 
$61440 \times 17784$ matrix ($8.1$~GB RAM) generated from video sequences from \url{http://perception.i2r.a-star.edu.sg/bk\_model/bk\_index.html}.
Such partial SVDs are used in the proximity operator for the nuclear norm term that arises in robust PCA formulations of video background subtraction. 
The randomized factorization happens to be faster than the Lanczos-based SVD from the PROPACK software even with a single core, but more importantly, the randomized method scales better as the parallelism increases. The accuracies of the two methods are indistinguishable. 
\label{fig:svdTiming}
}
%\vspace{-6mm}
\end{figure}

Figure \eqref{fig:svdTiming} illustrates the numerical benefits of such randomization over the classical Lanczos method. Since the randomized routine can perform all the multiplications in blocks, it benefits significantly from parallelization. 

 % spmag_parallel
\section{The role of parallel and distributed computation}

Thanks to Moore's law of scaling silicon density, raw computational throughput and storage capacity have increased at exponential rates up until the mid 2000's, thereby giving convex optimization algorithms commensurate performance boosts. However, while Moore's law is expected to continue for years to come, transistor \emph{efficiencies} have plateaued. As dictated by Dennard's law, scaling silicon density now results in unprecedented levels of power consumption. To handle the massive computational and storage resources demanded by Big Data at reasonable power costs, we must hence increasingly rely on parallel and distributed computation. 

While first-order methods seem ideally suited for near-optimal performance speed-ups, two issues block us when using distributed and heterogeneous hardware:
\begin{itemize}
\item {\bf Communication}: Uneven or faulty communication links between computers and within local memory hierarchy can significantly reduce the overall numerical efficiency of first-order methods. Two approaches broadly address such drawbacks. First, we can specifically design algorithms that minimize communication. Second, we can eliminate a master vector $\xk$ and instead work with a local copy in each machine that each lead to a consensus $\xStar$ at convergence.
\item {\bf Synchronization}: To \emph{exactly} perform the computations in a distributed fashion, first-order methods must coordinate the activities of different computers whose numerical primitives depend on the same vector $\xk$ at each iteration. However, this procedure slows down even when  a single machine takes much longer than the others. To alleviate this quintessential synchronization problem,  \emph{asynchronous} algorithms allow updates using outdated versions of their parameters. 
\end{itemize}

In this section we describe several key developments related to first-order methods within this context. Due to lack of space,  we will gloss over many important issues that impact the practical performance of these methods, such as latency  and multi-hop communication schemes. 

%\vspace{-3mm}
\subsection{Embarrassingly parallel first-order methods}
First-order methods can significantly benefit from parallel computing. These computing systems are typified by uniform processing nodes that are in close proximity and have reliable communications. Indeed, the expression \emph{embarrassingly parallel} refers to an ideal scenario for parallelization where we split the job into independent calculations that can be simultaneously performed in a predictable fashion. 

In parallel computing, the formulation of the convex problem makes a great deal of difference. An important embarrassingly parallel example is the computation of the gradient vector when the objective naturally decomposes as in \eqref{eq:SG}. Here, we can process each $F_i$ with one of $m$ computers using only $\order(n/m)$ local computation. Each machine also stores data locally with the corresponding $\order(n/m)$-data samples since each $F_i$ directly corresponds to a data point. Each processor then communicates with  the central location to form the final gradient and achieve the ideal linear speed-up. 

%\vspace{-3mm}
\begin{algorithm}[H]
\begin{algorithmic}[1]
    \newcommand{\xik}{\iterate{x_{(i)}}{k}}
    \newcommand{\xikk}{\iterate{x_{(i)}}{k+1}}
    \STATE $\zkk = \frac{1}{n}\sum_{i=1}^n \prox_{\gamma F_i}( \xik )$
    \FOR{$i=1$ \TO $n$}
        \STATE $\xikk = 2\zkk - \zk + \xik - \prox_{\gamma F_i}(\xik)$
        \ENDFOR
\end{algorithmic}
\caption{Decomposition algorithm (aka, consensus ADMM) \cite{CombettesPesquet2008} to solve \eqref{eq:sum_gi}; $\gamma>0$, $\xzero_i=0$ for $i=1,\ldots,n$.
}
\label{alg:4}
\end{algorithm}
%\vspace{-3mm}

Beyond parallelizing the basic gradient method for smooth problems, an embarrassingly parallel \emph{distribute and gather} framework for non-smooth problems results from an artificial reformulation of \eqref{eq:SG} so that we can apply decomposition techniques, such as Algorithm \ref{alg:4}:
\begin{equation}\label{eq:sum_gi}
\min_{x,x_{(i)}: i=1,\ldots,n} \left\{ \frac{1}{n}\sum_{i=1}^n F_i(x_{(i)}): x_{(i)} = x, i=1,\ldots,n\right\}.
\end{equation}
\noindent Indeed, the decomposition idea above forms the basis of the massively parallel consensus ADMM algorithm, which provides an extremely scalable optimization framework for $n>2$.
See \cite{CombettesPesquet2008,Condat2011,boyd2011distributed} for convergence analysis and further variants that include additional linear operators. 

Fortunately, we have access to many computer programming models to put these ideas immediately into action. Software frameworks, such as MapReduce, Hadoop, Spark, Mahout, MADlib, SystemML, and BigInsights, and corresponding high-level languages such as Pig, Hive and Jaql, 
can govern the various optimization tasks in parallel while managing all communications and data transfers within the computing system, and seamlessly provide for redundancy and fault tolerance in communications.

%\vspace{-3mm}
\subsection{First-order methods with reduced or decentralized communications}
In large systems, communicating the gradient or its elements to a central location may create a communication bottleneck. In this setting, coordinate descent methods provide a principled approach to reduce communications. There is indeed substantial work on developing parallel versions of these methods, dating back to work on the Jacobi algorithm for solving linear systems. The basic idea is simply to apply several coordinate descent updates at the same time in parallel. The advantage of this strategy in terms of communication is that each processor only needs to communicate a single coordinate update, while it only needs to receive the updates from the coordinates that have changed. 

When the objective is decomposable, this is simply an embarrassingly parallel version of the serial algorithm. Furthermore, classical work shows that this strategy is convergent, although it may require a smaller step size than the serial variant. However, it does not necessarily lead to a speed increase  for non-separable functions. Recent work has sought to precisely characterize the conditions under which parallel coordinate descent methods still obtain a large speed-up \cite{richtarik2012parallel}. 

Surprisingly, we can also decentralize the communication requirements of gradient methods for decomposable objectives with only minor modifications  \cite{shiextra}. The resulting algorithm performs a modified gradient update to the average of the parameter vectors only among the neighbors it communicates with. This strategy in fact achieves similar convergence rates to the gradient method with central communications, where the rate degradation depends on the graph Laplacian of the underlying communication network. 

%\vspace{-3mm}
\subsection{Asynchronous first-order methods with decentralized communications}
The gradient and the decomposition methods above still require a global synchronization to handle decomposable problems such as \eqref{eq:SG}. For instance, the gradient algorithm computes the gradient \emph{exactly} with respect to one (or more) examples at $\xk$ and then synchronizes in sequence to update $\xkk$ in a standard implementation. In contrast, stochastic gradient algorithms that address \eqref{eq:SG}  only use a crude approximation of the gradient. Hence, we expect these algorithm to be robust to outdated information, which can happen in asynchronous settings. 

A variety of recent works have shown that this is indeed the case. We highlight the work \cite{niu2011hogwild}, which models a lock-free shared-memory system where stochastic gradient updates are independently performed by each processor.
While the lock-free stochastic gradient still keeps a global vector $x$, processors are free to update it without any heed to other processors and continue their standard motions using the cached $x$. Under certain conditions this asynchronous procedure preserves the convergence of stochastic gradient methods, and results in substantial speed-ups when many cores are available. The same memory lock-free model also applies to  stochastic parallel coordinate descent methods \cite{richtarik2012parallel}. Finally, first-order algorithms with randomization can be effective even in asynchronous and decentralized settings with the possibility of communication failures~\cite{agarwal2011distributed}.

 % spmag_conclusion
\section{Outlook for convex optimization}
Big data problems necessitate a fundamental overhaul of how we design convex optimization algorithms, and suggest unconventional computational choices. To solve increasingly larger convex optimization problems with relatively modest growth in computational resources, this article makes it clear that we must identify key structure-dependent algorithmic approximation trade-offs. 

Since the synchronization and communication constraints of the available hardware naturally dictates the choice of the algorithms, we expect that new approximation tools will continue to be discovered that ideally adapt convex algorithms to the heterogeneity of computational platforms. 
We also predict an increased utilization of composite models and the corresponding proximal mapping principles for non-smooth Big Data problems to cope with noise and other constraints. 
For example, the LASSO formulation in \eqref{eq: lasso} has estimation guarantees that are quantitatively
stronger than the guarantees of the LS estimator when the signal $x_0$ has at most $k$ non-zero entries and $\Phi$ obeys certain assumptions~\cite{wainwright2014structured}.
That is to say, in order to get more out of the same data, we must use composite models.
This also invites the question of whether we can use composite models to get the same information out but do it faster, an issue which has been discussed~\cite{chandrasekaran2013computational,bottou2007tradeoffs} but not yet had an impact in practice.

\section*{Acknowledgements}
\small
Volkan Cevher's work is supported in part by the European Commission under grants MIRG-268398 and ERC Future Proof and by the Swiss Science Foundation under grants SNF 200021-132548, SNF 200021-146750, and SNF CRSII2-147633.
During the preparation of the work, Stephen Becker was supported as a Goldstein Fellow at the IBM T.\ J.\ Watson research center, and Mark Schmidt was supported by the Natural Language Laboratory at Simon Fraser University.

\small{
\setstretch{1} 
%\pdfbookmark[1]{References}{refSection}
\ifdefined\isdraft
\bibliographystyle{icml2013}
\else
\bibliographystyle{IEEEtran}
\fi
\bibliography{biblio-merged,biblio-volkan}

}

\end{document}